\documentclass[aos]{imsart}

\RequirePackage{amsthm,amsmath,amsfonts,amssymb}
\RequirePackage[numbers]{natbib}
\RequirePackage[colorlinks,citecolor=blue,urlcolor=blue]{hyperref}
\RequirePackage{graphicx}

\startlocaldefs
\theoremstyle{plain}

\newtheorem{remark}{Remark}
\newtheorem{corollary}{Corollary}[section]

\newtheorem{theorem}{Theorem}[section]
\newtheorem{lemma}{Lemma}[section]
\theoremstyle{remark}

\DeclareMathOperator*{\argmax}{arg\,max}
\DeclareMathOperator*{\argmin}{arg\,min}
\endlocaldefs

\begin{document}

\begin{frontmatter}
\title{Behaviour of FWER in Normal Distributions}
\runtitle{Behaviour of FWER in Normal Distributions}

\begin{aug}
\author[]{\fnms{Monitirtha} \snm{Dey}\ead[label=e1]{monitirtha.d\_r@isical.ac.in}}

\address[]{Indian Statistical Institute, 
\printead{e1}}

\end{aug}

\begin{abstract}
Familywise error rate (FWER) has been a cornerstone in simultaneous inference for decades, and the classical Bonferroni method has been one of the most prominent frequentist approaches for controlling FWER. The present article studies the behavior of the FWER for Bonferroni procedure in a multiple testing problem. We establish upper bounds on FWER for Bonferroni method under the equicorrelated and general normal setups in non-asymptotic case. 
\end{abstract}
 
\begin{keyword}[class=MSC]
\kwd[Primary ]{62J15}
\kwd[; secondary ]{62F03}
\end{keyword}

\begin{keyword}
\kwd{Multiple testing under dependence}
\kwd{Familywise error rate}
\kwd{Bonferroni method}
\end{keyword}

\end{frontmatter}

\section{Introduction}

Large-scale simultaneous inference problems in various disciplines often analyze related variables simultaneously. For example, in genome-wide association studies, high-density SNP (single nucleotide polymorphism) markers used to analyze genetic diversity exhibit high correlation. In spatial data with close geographical locations, the test statistics corresponding to different hypotheses often get influenced by each other. Multistage clinical trials and functional magnetic resonance imaging studies also concern variables with complex and unknown dependence structures. However, most classical multiple testing procedures controlling false discovery rate (FDR) or familywise error rate (FWER) typically rely on independence or some form of weak dependence among the concerned variables. Ignoring the dependence among related variables can produce highly variable significance measures and bias due to the confounding of dependent noise and the signal of interest~\cite{r11}. Efron~\cite{r6} mentions the correlation penalty on the summary
statistics depends on the root mean square (RMS) of correlations.

For these reasons, the problem of capturing the association among observations and extending the existing methods under association has attracted considerable attention in recent times. Multiple testing procedures under dependence have been discussed by Blanchard and Roquain~\cite{r1}, Efron~\cite{r5}, Liu, Zhang and Page~\cite{r13}, Sun and Cai~\cite{r17}, among others. Efron~\cite{r7} contains an excellent review of the relevant literature.

In this paper, we focus on the FWER, a widely considered frequentist approach in multiple testing. This is defined as the probability of making at least one false rejection in a family of hypothesis-testing problems. Controlling FWER has been a traditional concern in many multiple testing problems. This  tradition is reflected in the books by Hochberg and Tamhane~\cite{r8}, Westfall and Young~\cite{r20}, and the review by Tamhane~\cite{r18}. The control of FWER at some target level $\alpha$ requires each of the individual hypothesis to be tested at lower levels, e.g. in the Bonferroni procedure $\alpha$ is divided by the number of tests considered.

We have considered the equicorrelated normal distribution with positive correlation $\rho$ at first. Das and Bhandari~\cite{r3} have found that under this setup, FWER($\rho$) is a convex function of $\rho$ as the number of hypotheses grows to infinity. Consequently, they show that the FWER of the Bonferroni procedure is bounded by $\alpha(1-\rho)$, $\alpha$ being the desired level. Dey and Bhandari~\cite{r4} have shown that the Bonferroni FWER($\rho$) tends to zero asymptotically for any positive $\rho$. These works explicate the fact that Bonferroni’s procedure becomes very conservative for large-scale multiple testing problems under correlated setups. However, there is very little literature which elucidates the magnitude of the conservativeness of Bonferroni’s method in a dependent setup with small or moderate dimensions. 
In this paper, we address this research gap in a unified manner by establishing upper bounds on FWER of the Bonferroni method in the equicorrelated and arbitrarily correlated non-asymptotic setups.

Order statistics for exchangeable normal random variables have applications in biometrics~\cite{r16},~\cite{r19}. Also, the maximum of exchangeable normal random vector can be used to model the lifetime of parallel systems conveniently. The non-asymptotic bounds on FWER proposed in this work provide lower bounds on the c.d.f. of the failure time of the parallel systems. Loperfido~\cite{r14} has shown that the maximum of $n$ observations from equicorrelated normal distribution follows ($n-1$) dimensional skew normal distribution. Although the c.d.f. of multivariate skew normal distribution is very difficult to tackle, non-asymptotic bounds on the c.d.f. may be obtained along similar lines as in this article. 

This paper is organized as follows. In Section 2, we set up the framework and introduce the necessary notation. Section 3 contains theoretical results about the bounds on FWER in equicorrelated normal setup. Section 4 extends these results to arbitrarily dependent setups while Section 5 presents simulation findings. We conclude and discuss potential extensions of this work in Section 6. Proofs of the results are presented in the Appendix.

\section{Preliminaries} We consider a \textit{Gaussian sequence model}:
$$X_{i} {\sim} \mathcal{N}(\mu_{i},1),  \quad i=1, \ldots,n,$$
where $X_{i}$’s are independent and we are interested in the $n$ null hypotheses $H_{0i}:\mu_{i}=0$. The global null $ H_{0}=\bigcap_{i=1}^{n} H_{0 i}$ asserts that all means $\mu_i = 0$ for each $i$, while under the alternative, some $\mu_{i}$ is non-zero. We have considered one sided tests (i.e, $H_{0 i}$ is rejected for large values of $X_{i}$ (say $X_{i}>c$ for some cut-off $c$)). The most natural measure of type-I error in multiple testing is FWER, which is the probability of erroneously rejecting at least one true null hypothesis where the probability is computed under $H_{0}$. For common cut-off procedures (i.e procedures which use same cut-off for each of the $n$ hypotheses), this happens if $X_{i}>c$ for some $i$. So,
\begin{align*}FWER =\mathbb{P}\left(X_{i}>c\right. \text{for some} \left.i \mid H_{0}\right)
=\mathbb{P}_{H_{0}}\bigg(\bigcup_{i=1}^{n}\{X_i > c\}\bigg).\end{align*}

\noindent We have considered equicorrelated setup at first (in Section $3$) i.e, 
$$\operatorname{\mathbb{C}orr}\left(X_{i}, X_{j}\right)=\rho  \quad \forall i \neq j  \quad (\rho \geq 0).$$ 

\noindent In Section 4, we have dealt with general correlated setup:
$$\operatorname{\mathbb{C}orr}\left(X_{i}, X_{j}\right)=\rho_{ij}  \quad \forall i \neq j  \quad (\rho_{ij} \geq 0).$$

We consider FWER of Bonferroni’s method which, in the one-sided setting, rejects $H_{0i}$ if $X_{i}>\Phi ^{-1}(1-\alpha/n) (=c, \text{say})$ where $\alpha \in (0,1)$ is the desired level of significance and $\Phi$ is the c.d.f. of standard normal distribution. 

We shall denote the FWER of Bonferroni’s method under the equicorrelated normal setup with correlation $\rho$ and under a correlated normal setup with correlation matrix $R$ by $FWER(n , \alpha, \rho)$ and $FWER(n , \alpha, R)$, respectively. 

Das and Bhandari~\cite{r3} consider the equicorrelated framework and establish the following:
\begin{theorem}\label{thm2.1}
Suppose each $H_{0i}$ is being tested at size $\alpha_{n}$. If $\displaystyle \lim _{n \rightarrow \infty} n \alpha_{n}=\alpha \in(0,1)$ then, $FWER$ asymptotically is a convex function in $\rho \in [0,1]$.
\end{theorem}
For Bonferroni’s procedure, $\alpha_n = \alpha/n$ and thus Theorem \ref{thm2.1} also applies for Bonferroni’s method. Moreover, Theorem \ref{thm2.1} results in the following corollary.
\begin{corollary}\label{c2.1}
Given any $\alpha \in (0,1)$ and $\rho \in [0,1]$, $FWER(n, \alpha, \rho)$ is asymptotically bounded by $\alpha(1-\rho)$.
\end{corollary}
Corollary \ref{c2.1} shows that Bonferroni procedure controls FWER at a much smaller level than $\alpha$, when there is a large number of hypotheses. Dey and Bhandari~\cite{r4} prove a much stronger result than Corollary \ref{c2.1}. 
\begin{theorem}\label{thm2.2}
Given any $\alpha \in (0,1)$ and $\rho \in (0,1]$, $\lim_{n \to \infty} FWER(n, \alpha, \rho) = 0.$
\end{theorem}
They extend Theorem \ref{thm2.2} to arbitrarily correlated normal setups.
\begin{theorem}\label{thm2.3}
Let $\mathbf{\Sigma}_{n}$ be the correlation matrix of $X_{1}, \ldots, X_{n}$ with $(i,j)$’th entry $\rho_{ij}$ such that $\liminf \rho_{ij}=\delta>0$. Then, for any $\alpha \in (0,1)$,
$$\lim_{n \to \infty}FWER(n,\alpha,\mathbf{\Sigma}_{n}) = 0.$$
\end{theorem}

Theorem \ref{thm2.3} highlights the fundamental problem of using Bonferroni procedure in a multiple testing problem. We shall call the setup with very large number of hypotheses an asymptotic setup while setups with small or moderate number of hypotheses will be referred to as non-asymptotic setups. We summarize known and new results regarding behaviour of FWER of Bonferroni's procedure under various dependent normal setups in Table \ref{results}.
\begin{table*}
\caption{Results on Bonferroni FWER}
\label{results}
\begin{center}
\begin{tabular}{@{}lc@{}}
\hline
\textbf{Dependent Setup} & \textbf{Results on FWER}\\
 \hline
 Equicorrelated Asymptotic & Corollary \ref{c2.1}~\cite{r3}, Theorem \ref{thm2.2}~\cite{r4}  \\  
 General Asymptotic & Theorem \ref{thm2.3}~\cite{r4} \\ 
 Equicorrelated Non-asymptotic & Theorem \ref{thm3.1},\ref{thm3.3},\ref{thm3.4},\ref{thm3.5},\ref{thm3.6} and Corollary \ref{c3.2} \\ 
 General Non-asymptotic & Theorem \ref{thm4.1},\ref{thm4.2},\ref{thm4.3},\ref{thm4.4} and Corollary \ref{c4.2}  \\ \hline
\end{tabular}
 \end{center}
\end{table*}

\section{Bounds on FWER in General Non-asymptotic Setup}
In equicorrelated setups with small and moderate dimensions, the $\alpha(1-\rho)$ bound fails (we shall see this in detail in Section 5). We need large number of hypotheses, e.g 100 million to get values of FWER close to zero. Hence, establishing upper bounds on FWER in small and moderate number of hypotheses become relevant. The following result will be crucial towards this.
\begin{theorem}\label{thm3.1}
Under the equicorrelated normal set-up,
$$\text{FWER}(n, \alpha, \rho) \leq \alpha - \dfrac{n-1}{n}\cdot \dfrac{{\alpha}^2}{n}- \dfrac{n-1}{2\pi}\int_{0}^{\rho}\dfrac{1}{\sqrt{1-z^2}}e^{\frac{-{\Phi ^{-1}(1-\frac{\alpha}{n})}^2}{1+z}}dz.$$
\end{theorem}
\noindent It is noteworthy that this bound holds for any choice of $(n, \alpha)$ and any $\rho \geq 0$.

\begin{corollary}\label{3.1}
Under the equicorrelated normal set-up, if $\rho \leq \alpha/n$,
$\text{FWER}(n, \alpha, \rho) \leq \alpha - \dfrac{n-1}{n}\cdot \alpha\rho.$
\end{corollary}

Hence, throughout this work, we assume that $\rho \geq \alpha/n$. We observe that the bound mentioned in Theorem \ref{thm3.1} involves a definite integral which is very difficult to evaluate analytically. As we are interested in obtaining upper bounds for FWER, it is enough if we can find a lower bound to the integral. Towards this, we show the following theorem which will be crucial to obtain a lower bound to the integral mentioned in Theorem \ref{thm3.1}.
\begin{theorem}\label{thm3.2}
 Suppose $(X,Y)$ follows a bivariate normal distribution with parameters $(0,0,1,1,\rho)$ with $\rho \geq 0$. Then, for all $x \geq 2$,
$$\mathbb{P}(X \leq x, Y \leq x) \geq [\Phi(x)]^2+\dfrac{1}{2\pi}\cdot \sin^{-1}\rho\cdot e^{-\frac{x^2}{1+\frac{\rho}{2}}}.$$
\end{theorem}
Theorem \ref{thm3.2} can be used to establish the following corollary.
\begin{corollary}\label{c3.2}
Under the equicorrelated normal set-up, if $x=\Phi ^{-1}(1-\frac{\alpha}{n}) \geq 2$, 
$$FWER(n, \alpha, \rho) \leq \alpha - \dfrac{n-1}{n}\cdot \dfrac{{\alpha}^2}{n}- \dfrac{n-1}{2\pi}\cdot \sin^{-1}\rho\cdot e^{-\frac{x^2}{1+\frac{\rho}{2}}}.$$
\end{corollary}
We shall write $x$ for $\Phi ^{-1}(1-\alpha/n)$ from now on. We observe from simulation study that, the upper bound $\alpha(1-\rho)$ given by Corollary \ref{c2.1} holds for any nonnegative value of $\rho$ when $n \geq 10000$ and $\alpha \geq 0.01$. When $n = 10000$ and $\alpha =0.01$, we have $x=4.42$. This, along with the findings from our simulations suggest that the bound holds for $x \geq 4.42$. Therefore, here we restrict ourselves to the case $x \leq 4.42$. 

We also observe that, when $\rho \geq 0.5$, the bound $\alpha(1-\rho)$ works when $n \geq 900$ and $\alpha \geq 0.01$. When $n = 900$ and $\alpha =0.01$, we have $x=4.23$. This, along with the findings from our simulations suggest that, when $\rho \geq .5$, the bound works for $x \geq 4.23$. Therefore, when $\rho \geq .5$, we restrict ourselves to the case $x \leq 4.23$.

We shall also assume $\rho \geq 0.01$ from now on. We shall derive upper bounds on $FWER(n, \alpha, \rho)$ in each of the following four cases separately:
\begin{longlist}
    \item[Case 1.] $4.23 \geq x \geq 2, \rho \geq .5$
    \item[Case 2.] $4.42 \geq x \geq 2, .01 \leq \rho < .5$
    \item[Case 3.] $x \leq 2, \rho \geq .5$
    \item[Case 4.] $x \leq 2, \rho < .5$
\end{longlist}

\noindent \underline{Case 1. $4.23 \geq x \geq 2, \rho \geq .5$}
\begin{theorem}\label{thm3.3}
Let $4.23 \geq x \geq 2$ and $\rho \geq .5$. Then, 
$$\forall x \in [x_{l},x_{l+1}], \quad \quad \text{FWER}(n, \alpha, \rho) \leq \alpha - \dfrac{n-1}{n}\cdot \dfrac{{\alpha}^2}{n} -\frac{n-1}{n}\cdot \frac{\alpha \rho}{6} \cdot C_{x_{l}}$$
where $x_{l}$'s and $C_{x_{l}}$'s are as follows:
\begin{center}
\begin{tabular}{ |c|c|c|c|c|c|c| c|} 
 \hline
 $l$ & 1 & 2 & 3 & 4 & 5 & 6 & 7 \\ \hline
 $x_l$ & 2 & 2.56 & 3.06 & 3.33 & 3.71 & 3.93 & 4.23\\ \hline
 $C_{x_{l}}$ & 1 & $\frac{1}{2}$ & $\frac{1}{\pi}$ & $\frac{1}{2\pi}$ & $\frac{1}{\pi^{2}}$ & $\frac{1}{6\pi}$ & -\\ \hline
\end{tabular}
\end{center}
\end{theorem}
The proof follows from Corollary \ref{c3.2} and is included in the Appendix.
\vspace{2mm}

\noindent \underline{Case 2. $4.42 \geq x \geq 2, .01 \leq \rho < .5$}
\begin{theorem}\label{thm3.4}
 Let $4.42 \geq x \geq 2$ and $.01 \leq \rho < .5$. Let $I_{1}=[\frac{1}{3},.5)$, $I_{2}=[\frac{1}{2\pi},\frac{1}{3})$ and $I_{3}=[0.01,\frac{1}{2\pi})$. Then, for $\rho \in I_{i}$ with $i=1, 2, 3$ and for $x \in [x_m(i), x_{m+1}(i)]$,
$$\text{FWER}(n, \alpha, \rho) \leq \alpha - \dfrac{n-1}{n}\cdot \dfrac{{\alpha}^2}{n} -\frac{n-1}{n}\cdot \frac{\alpha \rho}{2\pi} \cdot D_{x_{m}}$$
where $x_{m}$'s and $D_{x_{m}}$'s are as follows:
\begin{center}
\begin{tabular}{ |c|c|c|c|c|c|c| c|c|c| c|} 
 \hline
 $m$ & 1 & 2 & 3 & 4 & 5 & 6 & 7 & 8 & 9 & \\ \hline
 $D_{x_{m}}$ & 1 & $\frac{1}{2}$ & $\frac{1}{\pi}$ & $\frac{1}{2\pi}$ & $\frac{1}{\pi^{2}}$ & $\frac{1}{\pi^{3}}$ & $\frac{1}{\pi^{4}}$ & $\frac{1}{4\pi^{4}}$ & $\frac{1}{16\pi^{4}}$ & \\ \hline
 $x_m(1)$ & 2 & 2.3 & 2.76 & 3 & 3.36 & 3.56 & 4 & 4.42 & &\\ \hline
 $x_m(2)$ &  & 2 & 2.49 & 2.72 & 3.04 & 3.23 & 3.66 & 4.03 &4.42 & \\ \hline
 $x_m(3)$ &  & 2 & 2.28 & 2.5 & 2.8 & 2.97 & 3.37 & 3.72 & 4.1 &4.42\\ \hline
\end{tabular}
\end{center}
\end{theorem}
The proof of this theorem is exactly similar to that of the previous theorem and hence omitted.

\vspace{3mm}
\noindent \underline{Case 3. $x \leq 2, \rho \geq .5$}
\begin{theorem}\label{thm3.5}
Let $x \leq 2$ and $\rho \geq .5$. Then, 
$$\text{FWER}(n, \alpha, \rho) \leq \alpha - \dfrac{n-1}{n}\cdot \dfrac{{\alpha}^2}{n} -\frac{n-1}{n}\cdot \frac{\alpha \rho}{6}.$$
\end{theorem}

\vspace{3mm}
\noindent \underline{Case 4. $x \leq 2, \rho < .5$}
\begin{theorem}\label{thm3.6}
Let $x =\Phi ^{-1}(1-\frac{\alpha}{n}) \leq 2$. Then, 
$$\text{FWER}(n, \alpha, \rho) \leq \alpha - \dfrac{n-1}{n}\cdot \dfrac{{\alpha}^2}{n} -\frac{n-1}{n}\cdot \frac{2 \alpha \rho}{5\pi}.$$
\end{theorem}
It is mention-worthy that Theorem \ref{thm3.6} is valid for any non-negative $\rho$. 


\section{Bounds on FWER in General Non-asymptotic Setup}
We have considered an equicorrelated dependence structure so far. However, problems involving variables with more general dependence structure need to be tackled with more general correlation matrices. Hence, the study of the behavior of FWER in arbitrarily correlated normal setups becomes crucial. Towards this, we consider the same \textit{Gaussian sequence model} as in Section 2, but now we assume $\operatorname{\mathbb{C}orr}\left(X_{i}, X_{j}\right)=\rho_{ij}$ for $i \neq j$ with $\rho_{ij} \geq 0$. Let $R$ be the correlation matrix of $X_{1}, \ldots, X_{n}$ and $FWER(n, \alpha,R)$ denote the FWER of Bonferroni’s method under this setup. So,
$$FWER(n, \alpha, R)=\mathbb{P}_{R}\bigg(\bigcup_{i=1}^n \{X_i > \Phi ^{-1}(1-\alpha/n) \}\mid H_{0}\bigg)=\mathbb{P}_{R}\left(\bigcup_{i=1}^{n} A_i\right)$$
where $A_i=\{X_i > \Phi ^{-1}(1-\frac{\alpha}{n}) | H_{0}\}$ for $i=1,\ldots, n$. 

In the equicorrelated setup, we use Kwerel's inequality (Lemma \ref{A.1}) to find an upper bound to FWER (see the Appendix):
$$\displaystyle \mathbb{P}\left(\bigcup_{i=1}^{n} A_i\right) \leq \sum_{i=1}^n \mathbb{P}(A_i) - \frac{2}{n} \sum_{1 \leq i < j \leq n} \mathbb{P}(A_i \cap A_j).$$
That approach can be used to obtain bounds on FWER in the arbitrarily correlated setup also. However, one observes that the above inequality gives equal importance to all the intersections. Therefore, it might be advantageous to use some other probability inequality which involves the intersections with higher probabilities only. We mention such an inequality below:
\begin{lemma}[Kounias~\cite{r10}]\label{4.1}
Let $A_1$, $A_2$, \ldots, $A_n$ be $n$ events. Then, $$\displaystyle \mathbb{P}\left(\bigcup_{i=1}^{n} A_i\right) \leq \sum_{i=1}^n \mathbb{P}(A_i) - \max_{1 \leq i \leq n} \sum_{j=1, j\neq i}^{n} \mathbb{P}(A_i \cap A_j).$$
\end{lemma}
Evidently Kounias's inequality is sharper than Kwerel's inequality and they are equivalent when $\mathbf{P}(A_{i}\cap A_{j})$ is same for all $i \neq j$. We are now in a position to state a generalization of Theorem \ref{thm3.1}:

\begin{theorem}\label{thm4.1}
 Consider the arbitrarily correlated normal set-up with covariance matrix $R$. Suppose $R$ has non-negative entries. Then, 
$$\text{FWER}(n, \alpha, R) \leq \alpha - \dfrac{n-1}{n}\cdot \dfrac{{\alpha}^2}{n}- \dfrac{1}{2\pi} \sum_{j=1, j\neq i^{*}}^{n} \int_{0}^{\rho_{i^{*}j}}\dfrac{1}{\sqrt{1-z^2}}e^{\frac{-{\Phi ^{-1}(1-\frac{\alpha}{n})}^2}{1+z}}dz$$
where $i^{*}=\displaystyle \argmax_{i} \sum_{j=1, j\neq i}^{n} \rho_{ij}$.
\end{theorem}

We observe that Theorem \ref{thm4.1} reduces to Theorem \ref{3.1} when $\rho_{ij}=\rho$ for all $i \neq j$.

\begin{corollary}\label{c4.1}
Consider the arbitrarily correlated normal setup with covariance matrix $R$. Suppose $R$ has non-negative entries. Let $i^{*}=\argmax_{i} \sum_{j=1, j\neq i}^{n} \rho_{ij}$ and $j_{*}=\argmin_{j}\rho_{i^{*}j}$. Then, if $\rho_{i^{*}j_{*}} \leq \frac{\alpha}{n}$,
$\text{FWER}(n, \alpha, R) \leq \alpha - \dfrac{n-1}{n}\cdot \alpha\rho_{i^{*}j_{*}}.$
\end{corollary}
Hence, we assume $\rho_{i^{*}j_{*}} > \frac{\alpha}{n}$ from now on. Suppose $\displaystyle \bar{\rho}_{i^{*}}=\frac{1}{n-1}\sum_{j=1, j\neq i^{*}}^{n} \rho_{i^{*}j}$. We have the following two generalizations of \ref{thm3.5} and \ref{thm3.6} respectively:

\begin{theorem}\label{thm4.2}
Consider the arbitrarily correlated normal setup with covariance matrix $R$. Suppose $R$ has non-negative entries. Let $\Phi ^{-1}(1-\frac{\alpha}{n}) \leq 2, \rho_{i^{*}j_{*}} \geq .5$. Then, 
$$\text{FWER}(n, \alpha, R) \leq \alpha - \dfrac{n-1}{n}\cdot \dfrac{{\alpha}^2}{n} -\frac{n-1}{n}\cdot \frac{\alpha \bar{\rho}_{i^{*}}}{6}.$$
\end{theorem}

\begin{theorem}\label{thm4.3}
Consider the arbitrarily correlated normal setup with covariance matrix $R$. Suppose $R$ has non-negative entries. Let $\Phi ^{-1}(1-\frac{\alpha}{n}) \leq 2$. Then, 
$$\text{FWER}(n, \alpha, R) \leq \alpha - \dfrac{n-1}{n}\cdot \dfrac{{\alpha}^2}{n} -\frac{n-1}{n}\cdot \frac{2\alpha \bar{\rho}_{i^{*}}}{5\pi}.$$
\end{theorem}

\noindent The proofs of these are exactly similar to those of Theorem \ref{thm3.5} and Theorem \ref{thm3.6} and hence omitted. 

In the proof of Theorem \ref{thm3.2}, we show that, for any $\rho \geq 0$,
$$\forall \hspace{.5mm} x \geq 2, \quad \quad \int_{0}^{\rho}\dfrac{1}{\sqrt{1-z^2}}e^{\frac{-x^2}{1+z}}dz \geq \sin^{-1} \rho \cdot e^{-\frac{x^2}{1+\frac{\rho}{2}}}.$$
This inequality leads to the following:
\begin{corollary}\label{c4.2}
Consider the arbitrarily correlated normal setup with covariance matrix $R$. Suppose $R$ has non-negative entries. Let $x = \Phi ^{-1}(1-\frac{\alpha}{n}) \geq 2$. Then, 
$$\text{FWER}(n, \alpha, R) \leq \alpha - \dfrac{n-1}{n}\cdot \dfrac{{\alpha}^2}{n}- \dfrac{1}{2\pi} \sum_{j=1, j\neq i^{*}}^{n} \sin^{-1}\rho_{i^{*}j}\cdot e^{-\frac{x^2}{1+\frac{\rho_{i^{*}j}}{2}}}$$
where $i^{*}=\argmax_{i} \sum_{j=1, j\neq i}^{n} \rho_{ij}$.
\end{corollary}

One can derive results similar to Theorem \ref{thm3.3} or Theorem \ref{thm3.4} using the above corollary by imposing certain conditions on the values of the correlations in the $i^{*}$-th row of $R$. For example, we have the following if we assume that $\rho_{i^{*}j_{*}} \geq .5$:
\begin{theorem}\label{thm4.4}
 Consider the arbitrarily correlated normal setup with covariance matrix $R$. Suppose $R$ has non-negative entries. Let $4.23 \geq x \geq 2$ and $\rho_{i^{*}j_{*}} \geq .5$. Then, 
$$\forall x \in [x_{l},x_{l+1}] \quad \quad \text{FWER}(n, \alpha, R) \leq \alpha - \dfrac{n-1}{n}\cdot \dfrac{{\alpha}^2}{n} -\frac{n-1}{n}\cdot \frac{\alpha \bar{\rho}_{i^{*}}}{6} \cdot C_{x_{l}}(n)$$
where $x_{l}$'s and $C_{x_{l}}(n)$'s are as follows:
\begin{center}
\begin{tabular}{ |c|c|c|c|c|c|c| c|} 
 \hline
 $l$ & 1 & 2 & 3 & 4 & 5 & 6 & 7 \\ \hline
 $x_l$ & 2 & 2.56 & 3.06 & 3.33 & 3.71 & 3.93 & 4.23\\ \hline
 $C_{x_{l}}(n)$ & 1 & $\frac{1}{2}$ & $\frac{1}{\pi}$ & $\frac{1}{2\pi}$ & $\frac{1}{\pi^{2}}$ & $\frac{1}{6\pi}$ & -\\ \hline
\end{tabular}
\end{center}
\end{theorem}
\noindent This can be established along the same lines of the proof of Theorem \ref{thm3.3}.

\section{Simulation Study}
The bound by Das and Bhandari~\cite{r3} provides a significant gain in power for Bonferroni
method for large number of hypotheses. However, for equicorrelated setups with small or moderate dimensions, their bound fails as mentioned earlier. We verify this through simulations. Our simulation scheme, for fixed $(n,\alpha)$ is as follows:
\begin{enumerate}
    \item For each $\rho \in \{0, .025, .050, .075, \ldots, 1\}$, we generate $10000$ $n$-variate equicorrelated multivariate normal observations (each with mean $0$ and variance $1$; common correlation coefficient being $\rho$).
    \item For each $\rho$, 
       \begin{itemize}
          \item in each of the 10000 replications, we note whether or not any of the generated $n$ components exceeds the cutoff ${\Phi}^{-1}(1-\alpha / n)$.
          \item the estimated FWER (for that $\rho$) is obtained accordingly from the 10000 replications.
       \end{itemize}
\end{enumerate}
We obtain the following plots after running these simulations for $(n,\alpha)=(100,.01)$ and $(500,.05)$ (the blue line represents the straight line $\alpha(1-\rho)$:
     
     
     \begin{figure}[h]
     \centering
     \begin{minipage}[b]{0.4\textwidth}
     \includegraphics[width=\textwidth]{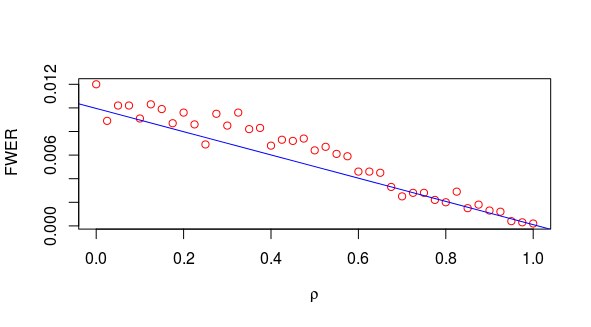}
     \end{minipage}
     \hfill
     \begin{minipage}[b]{0.4\textwidth}
     \includegraphics[width=\textwidth]{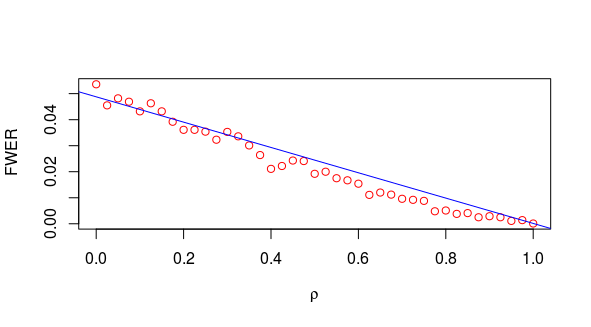}
     \end{minipage}
  \caption{FWER Plots for $(n,\alpha)$=(100,.01) and (500,.05)}
\end{figure}
\noindent We can see that the $\alpha(1-\rho)$ bound fails in these cases. Also, FWER is  not a convex function of $\rho$ in these cases. We present the simulation results for some choices of $(n,\alpha,\rho)$ along with our proposed bounds in Table~\ref{sphericcase}. It is mention worthy that in each case the estimated FWER is smaller than our proposed bounds.
 \begin{table}
\caption{Estimates of FWER($n,\alpha,\rho$)}
\label{sphericcase}
\begin{tabular}{@{}lrrrrrrc@{}}
\hline
 $(n,\alpha)$ & $x$ & Correlation ($\rho$) & 0.1 & 0.3 & 0.5 & 0.7 & 0.9 \\ \hline \hline
 $(10,0.3)$ & 1.8808 & $\widehat{\text{FWER}}(n, \alpha, \rho)$ & .2132 & .2053 & .1688 & .1242 & .0733 \\  
 &  & Bound  & .2156 & .2087 & .1965 & .1875 &.1785 \\ \hline
 $(100,0.05)$ & 3.2905 & $\widehat{\text{FWER}}(n, \alpha, \rho)$  & .0456 & .0355 & .0265 & .0153 & .0005 \\
  &  & Bound & .0475 & .0474 & .0462 & .0457 & .0452  \\ \hline
 $(500,0.05)$ & 3.7190 & $\widehat{\text{FWER}}(n, \alpha, \rho)$  & .0451 & .0319 & .0198 & .0081 & .0028  \\
  &  & Bound & .0475 & .0475 & .0471 & .0469 & .0467  \\ \hline
\end{tabular}
\end{table}

One can see that our bounds give good results for small values of equicorrelation $\rho$ and tend to become weak for large values of $\rho$. This is in contrast to the method of~\cite{r3} whose bound works in the large $\rho$ case. Therefore, in a way, our bounds and the $\alpha(1-\rho)$ bound are complementary to each other in depicting the behaviour of FWER in equicorrelated normal setups.

\section{Concluding Remarks}
This work is probably the first attempt in the context of finding the effect of correlation on FWER for small and moderate number of hypotheses.

The proofs of our results heavily use the fact that FWER can be regarded as $\mathbb{P}(\cup_{i=1}^{n} A_{i})$ for suitably defined events $A_i$, $1 \leq i \leq n$. Accurate computation of this probability is difficult because, in practice (as in multiple testing), the complete dependence between the events $(A_1,\ldots,A_n)$ is often unknown or unavailable (in our case $\rho$ is unknown and we have only some idea about $\rho$), unless the events $A_j$ are independent of each other. The available information is often the marginal probabilities and joint probabilities up to level $m(m<<n)$. In these situations, it is desirable to compute a lower or upper bound using only a limited amount of information. Our results utilize only the individual probabilities and the probabilities of pairwise intersections.

There are several interesting generalizations in various directions. One direction is to relax the multivariate normality assumption. Another direction is to use more general error rate criteria instead of FWER which control false rejections less severely, but in doing so are better able to detect false null hypotheses. In many areas, e.g microarray data analysis, the number of hypotheses under consideration is quite
large. Control of the FWER in those cases is so stringent that departures from the null hypothesis have little chance of being detected. Consequently, alternative measures of error control have been proposed in the literature. One such measure is $k$-FWER proposed by Lehmann and Romano~\cite{r12} which is the probability of rejecting at least $k$ true null hypotheses in a simultaneous testing problem. Such an error rate with $k>1$ is appropriate when one is willing to tolerate one or more false rejections, provided the number of false rejections is controlled. It is interesting to obtain similar upper bounds for $k$-FWER under arbitrarily dependent setups with small or moderate dimensions. 

\begin{appendix}
\section*{}
\subsection{Proof of Theorem \ref{thm3.1}}
We need two lemmas to establish this theorem.
\begin{lemma}[Kwerel~\cite{r9}]\label{A.1}
 Let $A_1$, $A_2$, \ldots, $A_n$ be $n$ events. Let $\displaystyle S_1= \sum_{i=1}^n \mathbb{P}(A_i)$ and $\displaystyle S_2=\sum_{1 \leq i < j \leq n} \mathbb{P}(A_i \cap A_j)$. Then, $\displaystyle \mathbb{P}\left(\bigcup_{i=1}^n A_i\right) \leq S_1 -\dfrac{2}{n}S_2$.
\end{lemma}

This bound on the union of $n$ events is also called the Sobel-Uppuluri upper bound and is the optimal linear bound in $S_1$ and $S_2$~\cite{r2}. The second lemma is regarding the joint distribution function of a bivariate normal distribution:

\begin{lemma}[Monhor~\cite{r15}]\label{A.2}
 Suppose $(X,Y)$ follows a bivariate normal distribution with parameters $(0,0,1,1,\rho)$ with $\rho \geq 0$. Then, for all $x >0$, 
$$\mathbb{P}(X \leq x, Y \leq x)=[\Phi(x)]^2+\dfrac{1}{2\pi}\int_{0}^{\rho}\dfrac{1}{\sqrt{1-z^2}}e^{\frac{-x^2}{1+z}}dz.$$
\end{lemma}

For $i=1,\ldots, n$, we define the event $A_i=\{X_i > \Phi ^{-1}(1-\alpha/n) | H_{0}\}$. 
So, $\mathbb{P}(A_i)=\mathbb{P}_{H_{0}}\left[X_i > \Phi ^{-1}(1-\alpha/n)\right]=\alpha/n$. This gives $S_1=\sum_{i=1}^n \mathbb{P}(A_i)=n\cdot \alpha/n=\alpha$. Now,
\begin{align*}
    &\mathbb{P}(A_{i} \cap A_{j})\\
    &=1 - \mathbb{P}(A_{i}^{c} \cup A_{j}^{c})\\
    &=1 - \mathbb{P}(A_{i}^{c}) - \mathbb{P}(A_{j}^{c}) + \mathbb{P}(A_{i}^{c} \cap A_{j}^{c})\\
    &=1-(1-\alpha/n)-(1-\alpha/n)+\mathbb{P}_{H_{0}}\bigg(X_i \leq \Phi ^{-1}(1-\alpha/n), X_j \leq \Phi ^{-1}(1-\alpha/n)\bigg)\\
    &=\frac{2\alpha}{n}-1+(1-\alpha/n)^2 +\dfrac{1}{2\pi}\int_{0}^{\rho}\dfrac{1}{\sqrt{1-z^2}}e^{\frac{-{\Phi ^{-1}(1-\frac{\alpha}{n})}^2}{1+z}}dz \quad \text{(using Lemma \ref{A.2})}\\
    &=\frac{\alpha^2}{n^2}+\dfrac{1}{2\pi}\int_{0}^{\rho}\dfrac{1}{\sqrt{1-z^2}}e^{\frac{-{\Phi ^{-1}(1-\frac{\alpha}{n})}^2}{1+z}}dz
\end{align*}
This gives 
$$S_2= \binom{n}{2} \cdot \bigg[\frac{\alpha^2}{n^2}+\dfrac{1}{2\pi}\int_{0}^{\rho}\dfrac{1}{\sqrt{1-z^2}}e^{\frac{-{\Phi ^{-1}(1-\frac{\alpha}{n})}^2}{1+z}}dz\bigg].$$
The rest is obvious from Lemma \ref{A.1} once we observe $FWER(n,\alpha,\rho) = \mathbb{P}(\bigcup_{i=1}^{n} A_i)$.
\vspace{2mm}

\subsection{Proof of Theorem \ref{thm3.2}} We use two well-known inequalities to prove this theorem.

\begin{lemma}[Chebyshev Integral Inequality]\label{A.3}
Let $f$ and $g$ be two nonnegative integrable functions and synchronous on a bounded interval $[a,b]$, i.e
$$\forall x,y \in [a,b], \quad \quad [f(x)-f(y)]\cdot[g(x)-g(y)] \geq 0.$$
Then,
$$(b-a)\cdot \int_{a}^{b}f(x)g(x)dx \geq \int_{a}^{b}f(x)dx \cdot\int_{a}^{b}g(x)dx.$$
\end{lemma}

\begin{lemma}[Hermite-Hadamard Integral Inequality]\label{A.4}
Let $f: [a,b] \to \mathbb{R}$ be a convex function. Then,
$$\int_{a}^{b}f(x)dx \geq (b-a)\cdot f\bigg(\frac{a+b}{2}\bigg).$$
\end{lemma}

Suppose $(X,Y)$ follows a bivariate normal distribution with parameters $(0,0,1,1,\rho)$, $\rho \geq 0$. Then, from Lemma \ref{A.2},
$$\forall x >0 \quad \quad \mathbb{P}(X \leq x, Y \leq x)=[\Phi(x)]^2+\dfrac{1}{2\pi}\int_{0}^{\rho}\dfrac{1}{\sqrt{1-z^2}}e^{\frac{-x^2}{1+z}}dz.$$
It can be easily shown that the functions $\dfrac{1}{\sqrt{1-z^2}}$ and $e^{\frac{-x^2}{1+z}}$ have same monotony in $z \in [0,1]$, i.e are synchronous on $[0,1]$. Using lemma \ref{A.3}, we obtain
\begin{equation}\int_{0}^{\rho}\dfrac{1}{\sqrt{1-z^2}}e^{\frac{-x^2}{1+z}}dz \geq \frac{1}{\rho} \int_{0}^{\rho}\dfrac{1}{\sqrt{1-z^2}}dz \cdot \int_{0}^{\rho}e^{\frac{-x^2}{1+z}}dz = \frac{\sin^{-1} \rho}{\rho} \cdot \int_{0}^{\rho}e^{\frac{-x^2}{1+z}}dz. \label{*} \tag{*}\end{equation}

The function $e^{\frac{-x^2}{1+z}}$ is convex in $z$ if $z \leq \frac{x^2}{2}-1$. Now, $0 \leq z \leq \rho \leq 1$. So, $z \leq \frac{x^2}{2}-1$ holds if $x \geq 2$. Hence, $e^{\frac{-x^2}{1+z}}$ is convex in $z \in [0,1]$ for $x \geq 2$. Applying Lemma \ref{A.4} on this function, we get,
$$\forall x \geq 2, \quad \quad \int_{0}^{\rho}e^{\frac{-x^2}{1+z}}dz \geq \rho \cdot e^{-\frac{x^2}{1+\frac{\rho}{2}}}.$$
Combining this with \eqref{*}, we get, for $x \geq 2$,
$$\int_{0}^{\rho}\dfrac{1}{\sqrt{1-z^2}}e^{\frac{-x^2}{1+z}}dz \geq \sin^{-1} \rho \cdot e^{-\frac{x^2}{1+\frac{\rho}{2}}}.$$
The rest is obvious from Lemma \ref{A.2}. 

\begin{remark}
Monhor~\cite{r15} obtained the following inequality for positively correlated bivariate normal distribution function using Lemma \ref{A.2}.
$$\mathbb{P}(X \leq x, Y \leq x) \geq [\Phi(x)]^2+\dfrac{1}{2\pi}\cdot \sin^{-1}\rho\cdot e^{-x^2} \quad \forall x>0 .$$ Theorem \ref{thm3.2} provides a sharper inequality for $x \geq 2$.
\end{remark}

\subsection{Proof of Theorem \ref{thm3.3}} We have from Corollary \ref{c3.2}, for each $x \geq 2$,
\begin{align*}
    FWER(n, \alpha, \rho) &\leq \alpha - \dfrac{n-1}{n}\cdot \dfrac{{\alpha}^2}{n}- \dfrac{n-1}{2\pi}\cdot \sin^{-1}\rho\cdot e^{-\frac{x^2}{1+\frac{\rho}{2}}}\\
    &\leq \alpha - \dfrac{n-1}{n}\cdot \dfrac{{\alpha}^2}{n}- \dfrac{n-1}{2\pi}\cdot \frac{2\pi\rho}{6}\cdot e^{-\frac{x^2}{1+.25}} \quad (\text{since}\hspace{1mm} \rho \geq .5 \hspace{1mm}\text{implies}\hspace{1mm}\frac{\sin^{-1}\rho}{\rho} \geq \frac{\pi}{3}).
\end{align*}
Hence it is enough to show that 
$$\forall x \in [x_{l},x_{l+1}] \quad \quad e^{-\frac{x^2}{1.25}} \geq \frac{\alpha}{n}\cdot C_{x_{l}}(n).$$
Now, $\frac{\alpha}{n}=1-\Phi(x)$. Let, $M(x)=\frac{e^{-\frac{x^2}{1.25}}}{1-\Phi(x)}$. Using computational tools, we get that $\forall x \in [x_{l},x_{l+1}]$, $M(x) \geq C_{x_{l}}(n)$ and the proof is completed. 

\subsection{Proof of Theorem \ref{thm3.5}} We have, from Theorem \ref{thm3.1},
$$FWER(n, \alpha, \rho) \leq \alpha - \dfrac{n-1}{n}\cdot \dfrac{{\alpha}^2}{n}- \dfrac{n-1}{2\pi}\int_{0}^{\rho}\dfrac{1}{\sqrt{1-z^2}}e^{\frac{-x^2}{1+z}}dz.$$
Now,
\begin{align*}
    &\dfrac{n-1}{2\pi}\int_{0}^{\rho}\dfrac{1}{\sqrt{1-z^2}}e^{\frac{-x^2}{1+z}}dz\\
    \geq \hspace{2mm} & \dfrac{n-1}{2\pi} \cdot \frac{1}{\rho} \int_{0}^{\rho}\dfrac{1}{\sqrt{1-z^2}}dz \cdot \int_{0}^{\rho}e^{\frac{-x^2}{1+z}}dz \quad \text{(using Lemma \ref{A.2})}\\
    = \hspace{2mm} & \dfrac{n-1}{2\pi} \cdot \frac{\sin^{-1}\rho}{\rho} \cdot  \bigg[\int_{0}^{\rho/2}e^{\frac{-x^2}{1+z}}dz+\int_{\rho/2}^{\rho}e^{\frac{-x^2}{1+z}}dz\bigg]\\
    \geq \hspace{2mm} & \dfrac{n-1}{2\pi} \cdot \frac{\sin^{-1}\rho}{\rho} \cdot \frac{\rho}{2} \bigg[e^{-x^2}+e^{-\frac{x^{2}}{{1+\rho/2}}}\bigg] \quad \text{(since $e^{\frac{-x^2}{1+z}}$ is increasing in $z$)} \\
    = \hspace{2mm} & \frac{\sin^{-1}\rho}{2\pi}\cdot(n-1) \cdot  \bigg[\dfrac{e^{-x^2}+e^{-\frac{x^{2}}{{1+\rho/2}}}}{2}\bigg]\\
    \geq \hspace{2mm} & \frac{\sin^{-1}\rho}{2\pi}\cdot(n-1) \cdot  \bigg[\dfrac{e^{-x^2}+e^{-\frac{x^{2}}{{1+.25}}}}{2}\bigg] \quad \text{(since $\rho \geq .5$)}\\
    = \hspace{2mm} & \frac{\sin^{-1}\rho}{2\pi}\cdot(n-1) \cdot G(x) \quad \text{(suppose)}
\end{align*}
Now, we have $\frac{\sin^{-1}\rho}{2\pi} \geq \frac{\rho}{6}$ since $\rho \geq .5$. Also, $G(x) \geq 1-\Phi(x)=\frac{\alpha}{n}$ for $x \leq 2.2$. The rest follows from Theorem \ref{thm3.1}. 

\subsection{Proof of Theorem \ref{thm3.6}}
$\Phi ^{-1}(1-\frac{\alpha}{n}) \leq 2$ implies $\frac{\alpha}{n} \geq 1-\Phi(2)=0.02275$. Therefore, $\rho \geq 0.02275$. Now, along the same lines of the preceding proof, we have,
\begin{align*}
    &\dfrac{n-1}{2\pi}\int_{0}^{\rho}\dfrac{1}{\sqrt{1-z^2}}e^{\frac{-x^2}{1+z}}dz\\
    \geq \hspace{2mm} & \frac{\sin^{-1}\rho}{2\pi}\cdot(n-1) \cdot  \bigg[\dfrac{e^{-x^2}+e^{-\frac{x^{2}}{{1+\rho/2}}}}{2}\bigg]\\
    \geq \hspace{2mm} & \frac{\rho}{2\pi}\cdot(n-1) \cdot  \bigg[\dfrac{e^{-x^2}+e^{-\frac{x^{2}}{{1+.011375}}}}{2}\bigg] \quad \text{(since $\rho \geq .02275$)}\\
    = \hspace{2mm} & \frac{\rho}{2\pi}\cdot(n-1) \cdot H(x) \quad \text{(suppose)}
\end{align*}
Now, $H(x) \geq \frac{4}{5}(1-\Phi(x))=\frac{4\alpha}{5n}$ for $x \leq 2$. The rest is obvious from Theorem \ref{thm3.1}. 

\subsection{Proof of Theorem \ref{thm4.1}} We have $\mathbb{P}(A_i)=\frac{\alpha}{n}$ where $A_i=\{X_i > \Phi ^{-1}(1-\frac{\alpha}{n}) | H_{0}\}$, for $i=1,\ldots, n$. One can show, along the similar lines of the proof of \ref{thm3.1}, the following: 

$$\mathbb{P}_{R}(A_i \cap A_j)=\frac{\alpha^2}{n^2}+\dfrac{1}{2\pi}\int_{0}^{\rho_{ij}}\dfrac{1}{\sqrt{1-z^2}}e^{\frac{-{\Phi ^{-1}(1-\frac{\alpha}{n})}^2}{1+z}}dz \quad \quad \forall \hspace{.5mm}i \neq j.$$
Hence, $\mathbb{P}_{R}(A_i \cap A_j)$ is an increasing function of $\rho_{ij}$. Therefore, $$\displaystyle \argmax_{i} \sum_{j=1, j\neq i}^{n} \mathbb{P}_{R}(A_i \cap A_j)=\argmax_{i} \sum_{j=1, j\neq i}^{n} \rho_{ij}= i^{*} \quad \text{(say)}.$$
Hence, applying Lemma \ref{4.1}, we get
\begin{align*}
    FWER=\displaystyle \mathbb{P}_{R}\left(\bigcup_{i=1}^{n} A_i\right)& \leq \sum_{i=1}^n \mathbb{P}(A_i) - \max_{1 \leq i \leq n} \sum_{j=1, j\neq i}^{n} \mathbb{P}_{R}(A_i \cap A_j)\\
    &= \alpha - \dfrac{n-1}{n}\cdot \dfrac{{\alpha}^2}{n}- \dfrac{1}{2\pi} \sum_{j=1, j\neq i^{*}}^{n} \int_{0}^{\rho_{i^{*}j}}\dfrac{1}{\sqrt{1-z^2}}e^{\frac{-{\Phi ^{-1}(1-\frac{\alpha}{n})}^2}{1+z}}dz
\end{align*}
, completing the proof. 
\end{appendix}
\begin{acks}[Acknowledgments]
The author would like to thank Prof. Subir Kumar Bhandari, ISI Kolkata, for illuminating discussions on controlling error rates in multiple hypothesis testing.
\end{acks}


\newpage
\begin{thebibliography}{4}
\bibitem{r1} \textsc{Blanchard, G.} and {Roquain, É.} (2009). Adaptive false discovery rate control under independence and dependence. \textit{J. Mach. Learn. Res.} \textbf{10} 2837–2871. MR2579914

\bibitem{r2} \textsc{Chen J.T.} (2014). \textit{Multivariate Bonferroni-Type Inequalities. Theory and Applications}. CRC Press.

\bibitem{r3} \textsc{Das, N.} and \textsc{Bhandari, S.K.} (2021). Bound on FWER for correlated normal. \textit{Statist. Probab. Lett.} \textbf{168} 108943. 

\bibitem{r4} \textsc{Dey, M.} and \textsc{Bhandari, S.K.} (2021). FWER Goes to Zero for Correlated Normal. \textit{arXiv preprint arXiv:2110.05070}. 

\bibitem{r5} \textsc{Efron, B.} (2007). Correlation and large-scale simultaneous significance testing. \textit{J. Amer. Statist. Assoc.} \textbf{102} 93–103. 

\bibitem{r6} \textsc{Efron, B.} (2010). Correlated z-values and the accuracy of large-scale statistical estimates. \textit{J. Amer. Statist. Assoc.} \textbf{105} 1042–1055. 

\bibitem{r7} \textsc{Efron, B.} (2010). \textit{Large-Scale Inference. Institute of Mathematical Statistics (IMS) Monographs} \textbf{1.} Cambridge Univ. Press, Cambridge. 

\bibitem{r8} \textsc{Hochberg, Y.} and \textsc{Tamhane, A.C.} (1987). \textit{Multiple Comparison Procedures}. Wiley, New York.  

\bibitem{r9} \textsc{Kwerel, S.M.} (1975). Most stringent bounds on aggregated probabilities of partially specified dependent probability systems. \textit{J. Amer. Statist. Assoc.} \textbf{70} 472–479.

\bibitem{r10} \textsc{Kounias, EG.} (1968). Bounds for the Probability of a Union, with Applications. \textit{Ann. Math. Stat.} \textbf{39} 2154-2158.

\bibitem{r11} \textsc{Leek, J.T.} and \textsc{Storey, J.D.} (2008). A general framework for multiple testing dependence. \textit{Proc.
Natl. Acad. Sci. USA} \textbf{105} 18718-18723.

\bibitem{r12} \textsc{Lehmann, E.L.} and \textsc{Romano, J.P.} (2005). 
 Generalizations of the familywise error rate. \textit{Ann. Statist.} \textbf{33} 1138 – 1154. 

\bibitem{r13} \textsc{Liu, J.}, \textsc{Zhang, C.} and \textsc{Page, D.} (2016). Multiple testing under dependence via graphical models. \textit{Ann. Appl. Stat.} \textbf{10} 1699-1724. 

\bibitem{r14} \textsc{Loperfido, N.}, \textsc{Navarro, J.}, \textsc{Ruiz, J.M.}, and \textsc{Sandoval, C.J.} (2007). Some relationships between skew-normal distributions and order statistics from exchangeable normal random vectors. \textit{Comm. Statist. Theory Methods} \textbf{36} 1719–1733. 

\bibitem{r15} \textsc{Monhor, D.} (2013). Inequalities For Correlated Bivariate Normal Distribution Function. \textit{Probab. Engrg. Inform. Sci.} \textbf{27} 115-123.

\bibitem{r16} \textsc{Olkin, I.} and \textsc{Viana, M.} (1995). Correlation analysis of extreme observations from a multivariate normal distribution. \textit{J. Amer. Statist. Assoc.} \textbf{90} 1373–1379.

\bibitem{r17} \textsc{Sun, W.} and \textsc{Cai, T. T.} (2009). Large-scale multiple testing under dependence. \textit{J. R. Stat. Soc. Ser. B. Stat. Methodol.} \textbf{71} 393–424. 

\bibitem{r18} \textsc{Tamhane, A.C.} (1996). Multiple comparisons. In \textit{Handbook of Statistics} (S. Ghosh and C. R. Rao, eds.) \textbf{13} 587–629. North-Holland, Amsterdam.

\bibitem{r19} \textsc{Viana, M.} (1998). Linear combinations of ordered symmetric observations with application
to visual acuity. In \textit{Handbook of Statistics} (N. Balakrishnan and C. R. Rao, eds.) \textbf{17} 513–524. North-Holland, Amsterdam.

\bibitem{r20} \textsc{Westfall, P.H.} and \textsc{Young, S.S.} (1993). \textit{Resampling-Based Multiple Testing}. Wiley, New York.




\end{thebibliography}
\end{document}